\documentstyle[12pt,amscd]{amsart}
\newtheorem{thm}[equation]{Theorem}
\numberwithin{equation}{section}

\newtheorem{lem}[equation]{Lemma}

\newtheorem{diag}[equation]{Diagram}
\newtheorem{prop}[equation]{Proposition}
\newtheorem{tab}[equation]{Table}

\begin{document}
\raggedbottom \voffset=-.7truein \hoffset=0truein \vsize=8truein
\hsize=6truein \textheight=8truein \textwidth=6truein
\baselineskip=18truept
\def\l{\lambda}
\def\ss{\smallskip}
\def\la{\langle}
\def\ra{\rangle}
\def\on{\operatorname}
\def\a{\alpha}
\def\bz{{\Bbb Z}}
\def\im{\on{im}}
\def\ext{\on{Ext}}
\def\sq{\on{Sq}}
\def\eps{\epsilon}
\def\si{\sigma}
\def\tfrac{\textstyle\frac}
\def\w{\wedge}
\def\equ{\begin{equation}}
\def\b{\beta}
\def\G{\Gamma}
\def\g{\gamma}
\def\endeq{\end{equation}}
\def\A{{\cal A}}
\def\P{{\cal P}}
\def\cj{{\cal J}}
\def\zt{{\bold Z}_2}
\def\Hom{\on{Hom}}
\def\ker{\on{ker}}
\def\coker{\on{coker}}
\def\da{\downarrow}
\def\colim{\operatornamewithlimits{colim}}
\def\io{\iota}
\def\Om{\Omega}
\def\u{{\cal U}}
\def\e{{\cal E}}
\def\exp{\on{exp}}
\def\wbar{{\overline w}}
\def\xbar{{\overline x}}
\def\ybar{{\overline y}}
\def\zbar{{\overline z}}
\def\ebar{{\overline e}}
\def\nbar{{\overline n}}
\def\rbar{{\overline r}}
\def\kt{{\widetilde k}}
\def\ni{\noindent}
\def\coef{\on{coef}}
\def\den{\on{den}}
\def\ot{\otimes}
\def\ms{\medskip}
\def\dirlim{\operatornamewithlimits{dirlim}}
\def\lfl{\lfloor}
\def\rfl{\rfloor}
\def\ord{\on{ord}}
\def\gd{{\on{gd}}}
\def\rk{{{\on{rk}}_2}}
\def\bsp{\on{BSp}}
\def\bspt{{\widetilde\bsp}}
\def\but{{\widetilde{BU}}}
\def\N{{\Bbb N}}
\def\Z{{\Bbb Z}}
\def\Q{{\Bbb Q}}
\def\R{{\Bbb R}}
\def\C{{\Bbb C}}
\def\buo{B^{\text{u}}}
\def\boo{B^{\text{o}}}
\def\bg{\bigg}
\def\mo{\on{mod}}
\def\Remark{\noindent{\it  Remark}}
\title[Immersions of complex projective space]
{Some new immersion results for complex projective space}
\author{Donald M. Davis}
\address{Department of Mathematics, Lehigh University\\Bethlehem, PA 18015, USA}
\email{dmd1@@lehigh.edu}
\date{February 21, 2006}

\keywords{Immersions, complex projective space, obstruction theory}
\thanks {2000 {\it Mathematics Subject Classification}:
57R42, 55S35.}

\maketitle
\begin{abstract} We prove the following two new optimal immersion results for complex projective space.
 First, if $n\equiv3$ mod 8 but $n\not\equiv3$ mod 64, and
$\a(n)=7$, then
$CP^{n}$ can be
immersed in $\R^{4n-14}$. Second, if $n$ is even and $\a(n)=3$, then $CP^n$ can be immersed in $\R^{4n-4}$.
Here $\a(n)$ denotes the number of 1's in the binary expansion of $n$.
The first contradicts a result of Crabb, which said that such an immersion
does not exist, apparently due to an arithmetic mistake. We combine Crabb's method with that developed by the author and Mahowald.
 \end{abstract}

\section{Main theorems}\label{intro}
We prove the following two new optimal immersion results for $2n$-dimensional complex projective space $CP^n$.
\begin{thm}\label{mainthm} If $n\equiv3$ mod $8$ and $n\not\equiv3$ mod $64$, and $\a(n)=7$,
then 
$CP^n$ can be immersed in $\R^{4n-14}$.\end{thm}
\begin{thm}\label{thm2} If $n$ is even and $\a(n)=3$, then $CP^n$ can be immersed in $\R^{4n-4}$.\end{thm}
Here and throughout,
$\a(n)$ denotes the number of 1's in the binary expansion of $n$.
Theorem \ref{mainthm} contradicts a result of Crabb (\cite{Crabb}). In Section \ref{sec2},
we prove Theorem \ref{mainthm} by an adaptation of Crabb's argument, and point out what we believe to be
his mistake, apparently arithmetic.
We prove Theorem \ref{thm2} in Section \ref{thm2sec}.

We now summarize what we feel to be the status of the immersion question for $CP^n$. In addition to
incorporating the two new immersion results above, we list as unresolved one immersion result claimed by Crabb.
We will discuss our reason for doing so in Section \ref{sec4}. Despite our feeling that two of Crabb's many results
are flawed, we feel that his overall approach, combining K-theory with obstruction theory, is sound;
we have checked the details of his immersion results cited in Theorem \ref{nodd}. Now we begin our summary.

 There are three families of results that apply to all values of $\a(n)$. All known nonimmersion
results follow from the first two. 
\begin{thm} $(\cite{SchwS})$ \label{SSMthm}  $CP^n$ cannot
be immersed in $\R^{4n-2\a(n)+\eps}$, where
$$\eps=\begin{cases}0&\text{if $n$ is even and $\a(n)\equiv1\ (\mo\,4)$}\\
1&\text{if $n$ is even and $\a(n)\equiv2,3\ (\mo\,4)$}\\
-1&\text{otherwise}.\end{cases}$$
 \end{thm}

\begin{thm} $(\cite{SigS})$ If $CP^n$ immerses in $\R^{4n-2\a(n)}$, then
$$\nu(c_{\a(n)-1})=\nu(c_{\a(n)})<\nu(c_i)\text{\quad for all }i<\a(n)-1,$$
where $\sum c_it^i=((\log(1+t))/t)^{2n+1-\a(n)}.$\label{SigSthm}\end{thm}
\noindent Here and throughout, $\nu(-)$ denotes the
exponent of 2 dividing an integer. The specific results obtainable from Theorem
\ref{SigSthm} were determined for $\a(n)\le5$ in \cite{SigS} and for $\a(n)=6$ and 7
(with a mistake when $\a(n)=7$) in \cite{Crabb}. In Section \ref{sec5}, we derive these
for $\a(n)=8$; the results have been incorporated into Tables \ref{evtabl} and \ref{odtabl}.

For large values of $\a(n)$, the best immersion results are obtained in the following embedding theorem,
which relied on earlier results of Milgram.
\begin{thm} $(\cite{Muk})$ $CP^n$ can be embedded in $\R^{4n-\a(n)}$, and, if $n>1$ is odd, $CP^n$ can be embedded in
$\R^{4n-\a(n)-1}$.\end{thm}

For small $\a(n)$, better immersion results are obtained by \cite{DM} and \cite{Crabb}. Crabb did not consider
even values of $n$, and so, when $n$ is even, the immersions are from \cite{DM} and Theorem \ref{thm2}, and the nonimmersions from
\cite{SchwS} (stated here as \ref{SSMthm}) and from \cite{SigS} (stated here as \ref{SigSthm}). 
\begin{thm} If $n$ is even, then $CP^n$ immerses in $\R^{4n-d}$ and does not immerse in $\R^{4n-e}$, where $d$ and $e$
are given in Table \ref{evtabl}.\end{thm}
\begin{center}
\begin{minipage}{6.5in}
\begin{tab}\label{evtabl}{\bf Immersions and nonimmersions when $n$ is even}
\begin{center}
\begin{tabular}{cl|rr}
$\a(n)$&&$d$&$e$\\
\hline
$2$&&$2$&$3$\\
$3$&&$4$&$5$\\
$4$&$n\not\equiv10\ (\mo\,16)$&$7$&$8$\\
$4$&$n\equiv10\ (\mo\,16)$&$7$&$9$\\
$5$&&$9$&$10$\\
$6$&&$10$&$11$\\
$7$&&$11$&$13$\\
$8$&$n\not\equiv22\,(32)\text{ and }\not\equiv132\,(256)$&$14$&$16$\\
$8$&$n\equiv22\,(32)\text{ or }\equiv132\ (256)$&$14$&$17$\\
$>8$&&$14$&
\end{tabular}
\end{center}
\end{tab}
\end{minipage}
\end{center}
\noindent Thus, when $n$ is even, the only unresolved case for $\a(n)\le6$ occurs when $\a(n)=4$ and $n\equiv10$ mod 16.

We believe that the following tabulation of results and earliest proofs is accurate
when $n$ is odd. Note that the case discussed in Proposition \ref{unresolved} is the only unresolved case
when $n$ is odd and $\a(n)\le 7$.
\begin{thm} If $n$ is odd, then $CP^n$ immerses in $\R^{4n-d}$ and does not immerse in $\R^{4n-e}$, where $d$ and $e$ are given in Table \ref{odtabl}.\label{nodd}\end{thm}
\begin{center}
\begin{minipage}{6.5in}
\begin{tab}\label{odtabl}{\bf Immersions and nonimmersions when $n$ is odd}
\begin{center}
\begin{tabular}{cl|rc|rl}
$\a(n)$&condition&$d$&ref&$e$&ref\\
\hline
$2$&&$3$&\cite{James}&$4$&\cite{SigS}\\
\hline
$3$&$n\equiv1\ (\mo\,4)$&$5$&\cite{Ran}&$6$&\cite{SigS}\\
$3$&$n\equiv3\ (\mo\,4)$&$6$&\cite{Steer}&$7$&\cite{SchwS}\\
\hline
$4$&$n\not\equiv7\ (\mo\,8)$&$7$&\cite{DM}&$8$&\cite{SigS}\\
$4$&$n\equiv7\ (\mo\,8)$&$8$&\cite{Crabb}&$9$&\cite{SchwS}\\
\hline
$5$&$n\equiv1\ (\mo\,4)$&$9$&\cite{DM}&$11$&\cite{SchwS}\\
$5$&$n\equiv3\ (\mo\,8)$&$9$&\cite{DM}&$10$&\cite{SigS}\\
$5$&$n\equiv7\ (\mo\,8)$&$10$&\cite{Crabb}&$11$&\cite{SchwS}\\
\hline
$6$&$n\equiv1\ (\mo\,4)$&$11$&\cite{DM}&$12$&\cite{SigS}\\
$6$&$n\equiv3\ (\mo\,16)$&$11$&\cite{DM}&$12$&\cite{SigS}\\
$6$&$n\equiv11\ (\mo\,16)$&$12$&\cite{Crabb}&$13$&\cite{SchwS}\\
$6$&$n\equiv7\ (\mo\,8)$&$12$&\cite{Crabb}&$13$&\cite{SchwS}\\
\hline
$7$&$n\equiv1\ (\mo\,4)$&$13$&\cite{DM}&$14$&\cite{SigS}\\
$7$&$n\equiv3\ (\mo\,64)$&$13$&\cite{DM}&$14$&\cite{SigS}\\
$7$&$n\equiv3\,(8),\not\equiv3\ (64)$&$14$&Thm.\ref{mainthm}&$15$&\cite{SchwS}\\
$7$&$n\equiv7\ (\mo\,8)$&$14$&\cite{Crabb}&$15$&\cite{SchwS}\\
\hline
$8$&$n\not\equiv15\,(16)\text{ and }\not\equiv37\,(64)$&$15$&\cite{DM}&$16$&\cite{SigS}\\
$8$&$n\equiv15\,(16)\text{ or }\equiv37\,(64)$&$15$&\cite{DM}&$17$&\cite{SchwS}\\
\hline
$>8$&&$15$&\cite{DM}&&\\
\hline
\end{tabular}
\end{center}
\end{tab}
\end{minipage}
\end{center}
\section{Proof of Theorem \ref{mainthm}}\label{sec2}
In this section we prove Theorem \ref{mainthm} and describe what we believe was Crabb's mistake when
he asserted a nonimmersion in this situation.

Let $H^\C_n$ denote the Hopf bundle over $CP^n$. It is standard that the immersion of \ref{mainthm} is equivalent
to showing that the stable normal bundle
$-(n+1)H_n^\C$ is stably equivalent to a bundle of dimension $2n-14$. We let
$n=8p+3$ with $\a(p)=5$. In \cite[\S3, esp.~(3.2)]{Crabb}, Crabb showed that a necessary condition for the immersion is
that, if $\l(T)=(\sinh^{-1}(\sqrt T)/\sqrt T)^2$, and $(\l^{8p})_i$ denotes the coefficient of
$T^i$ in $\l(T)^{8p}$, then there exists an integer $e$ such that $e(\l^{8p})_3\equiv 64$ mod 128,
$e(\l^{8p})_2\equiv0$ mod 32, $e(\l^{8p})_1\equiv0$ mod 8, and $e(\l^{8p})_0\equiv0$ mod 2.
(For the reader wishing to compare with Crabb's notation, his $l=14$, $j=3$, and $k=4p+1$.)

Working mod $(T^4)$, we have 
$$\l(T)=1-\tfrac13T+\tfrac8{45}T^2-\tfrac4{35}T^3,$$
and
$$\l(T)^8=1-\tfrac83T+\tfrac{68}{15}T^2-\tfrac{1192}{189}T^3.$$
We then have
$$\l(T)^{8p}=1+u_18pT+u_24pT^2+u_38pT^3,$$
where each $u_i$ is an odd fraction.
The first and last of Crabb's necessary conditions stated above require
$ep\equiv8$ mod 16 and $e\equiv0$ mod 2. These (and the other conditions)
can be satisfied if and only if $p\not\equiv0$ mod 8. Crabb's Lemma 3.4 makes it 
clear that he believed that his conditions could not
 be satisfied in the cases which
we address here.

Now we prove that the immersion exists when $p\not\equiv0$ mod 8. We use modified Postnikov
towers (MPTs), as introduced in \cite{GM} and employed in many papers  such as
\cite{DM1}, \cite{DM}, and, more recently, \cite{Sin}. We consider the lifting question
\begin{equation}\label{lift?}\begin{CD}@.@.\bspt(16p-8)\\
@.@.@VqVV\\ CP^{8p+3}@>h>> HP^{4p+1} @>f>>\bsp\end{CD}\end{equation}
where $f$ classifies the stable bundle $-(4p+2)H^{Sp}_{4p+1}$ over the quaternionic projective
space. The space $\bspt(m)$ is the classifying space for symplectic vector bundles of real geometric dimension $m$.
It is the pullback of $BO(m)$ and $\bsp$ over $BO$. We let
$$\bspt(16p-8)=E_8\to E_7\to \cdots \to E_1\to \bsp$$
denote the MPT through dimension $16p+6$. In this range, the fiber of $q$ is the stable stunted
real projective space $P_{16p-8}=RP^\infty/RP^{16p-9}$, whose homotopy groups in this range are displayed in
\cite[Table 8.9]{Mem}. We reproduce them in Diagram \ref{chart}, indexed as $\pi_*(\Sigma P_{16p-8})$, which is their
dimensions as $k$-invariants in the MPT.

\begin{center}
\begin{minipage}{6.5in}
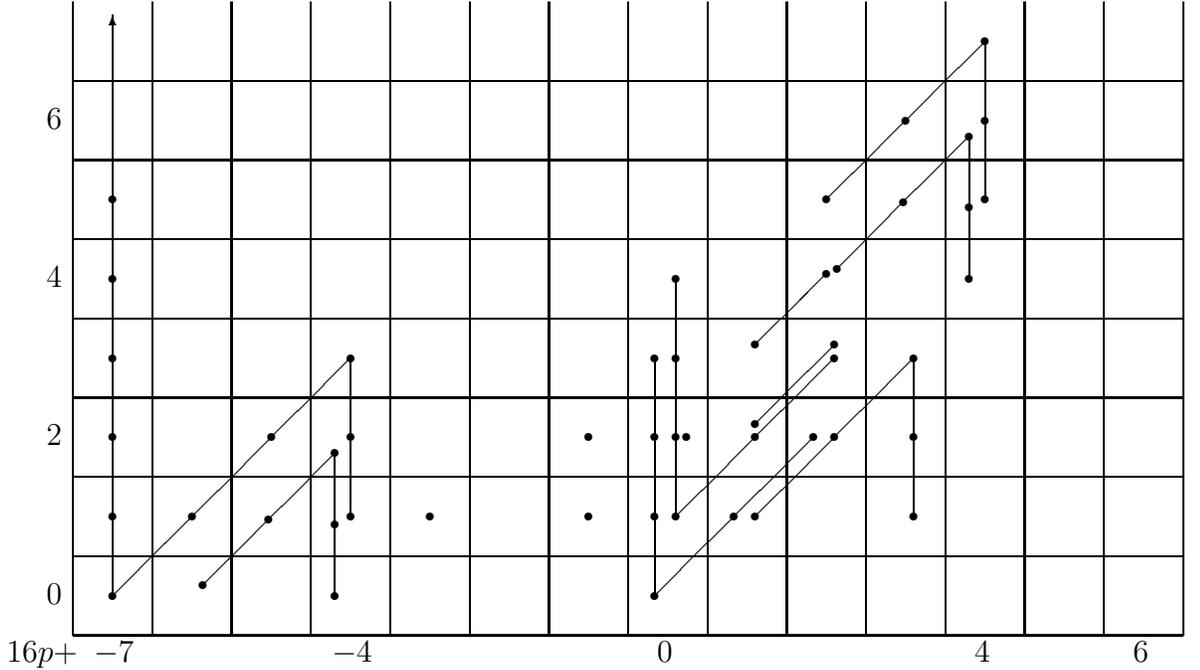
\begin{diag}\label{chart}{\bf Adams spectral sequence of $\Sigma P_{16p-8}$}
\begin{center}
\begin{picture}(445,260)
\def\mp{\multiput}
\def\elt{\circle*{3}}
\put(0,0){$16p+$}
\put(33,0){$-7$}
\put(123,0){$-4$}
\put(246,0){$0$}
\put(366,0){$4$}
\put(426,0){$6$}
\put(15,22){$0$}
\put(15,82){$2$}
\put(15,142){$4$}
\put(15,202){$6$}
\mp(25,10)(0,30){8}{\line(1,0){420}}
\mp(25,10)(30,0){15}{\line(0,1){240}}
\mp(40,25)(0,30){6}{\elt}
\put(40,25){\vector(0,1){220}}
\mp(70,55)(30,30){3}{\elt}
\put(40,25){\line(1,1){90}}
\mp(130,55)(0,30){2}{\elt}
\mp(130,55)(240,120){2}{\line(0,1){60}}
\mp(370,175)(0,30){3}{\elt}
\mp(74,29)(25,25){3}{\elt}
\mp(124,25)(0,27){2}{\elt}
\mp(74,29)(240,120){2}{\line(1,1){50}}
\mp(124,25)(240,120){2}{\line(0,1){54}}
\mp(314,149)(25,25){3}{\elt}
\mp(364,145)(0,27){2}{\elt}
\mp(310,175)(30,30){2}{\elt}
\put(310,175){\line(1,1){60}}
\mp(245,25)(0,30){4}{\elt}
\mp(245,25)(8,30){2}{\line(0,1){90}}
\mp(253,55)(0,30){4}{\elt}
\mp(245,25)(8,30){2}{\line(1,1){60}}
\mp(275,55)(30,30){2}{\elt}
\mp(283,85)(30,30){2}{\elt}
\mp(283,55)(30,30){3}{\elt}
\mp(343,55)(0,30){2}{\elt}
\put(283,55){\line(1,1){60}}
\put(343,55){\line(0,1){60}}
\mp(283,90)(30,30){2}{\elt}
\put(283,90){\line(1,1){30}}
\put(283,120){\line(1,1){27}}
\mp(283,120)(27,27){2}{\elt}
\mp(160,55)(60,0){2}{\elt}
\mp(220,85)(37,0){2}{\elt}
\end{picture}
\end{center}
\end{diag}
\end{minipage}
\end{center}

\medskip
The obstructions for lifting from $E_i$ to $E_{i+1}$ are $k$-invariants in $H^j(E_i)$
corresponding to dots in position $(j,i)$ of the diagram. All cohomology groups have coefficients in $\Z/2$.
The bulk of our work will be in
proving the following result, which states that $f$ lifts to the fifth stage of the MPT.
\begin{prop}\label{E5} In $(\ref{lift?})$, $f$ factors through a map $HP^{4p+1}@>f_5>>E_5$.\end{prop}

Before giving the proof of Proposition \ref{E5}, we use it to complete the proof of Theorem \ref{mainthm}.
Let $\ell=f_5\circ h:CP^{8p+3}\to E_5$.
To get $CP^{8p+1}$ to lift to $\bspt(16p-8)$, we need only show that $\ell^*(k_2)=0$, where $k_2\in H^{16p+2}(E_5)$
corresponds to the dot in position $(16p+2,5)$. The diagonal line emanating from this dot suggests, and the
computation of the MPT proves, that there is a relation in $H^*(E_5)$ of the form $\sq^2k_2+ak_{-7}=0$,
where $k_{-7}\in H^{16p-7}(E_5)$ corresponds to the dot at height 5 in the initial tower, and $a$ is a combination of
Steenrod operations and Stiefel-Whitney classes acting on $k_{-7}$. Therefore, since $H^{16p-7}(CP^{8p+3})=0$,
we must have $\sq^2(\ell^*(k_2))=0$ in $H^*(CP^{8p+3})$. Since $\sq^2$ acts injectively on $H^{16p+2}(CP^{8p+3})$,
this implies $\ell^*(k_2)=0$, and hence $CP^{8p+1}$ lifts to $\bspt(16p-8)$. 

By \cite[Prop.3.2]{Crabb}, the $KO$-theoretic obstruction to extending this lifting over $CP^{8p+2}$ is given by the conditions
on $\l^{8p}$ described above, which are satisfied under our hypotheses, and hence this $KO$-theoretic
obstruction is 0.
The total obstruction to this extension lies in $\pi_{16p+4}(\Sigma P_{16p-8})\approx\Z/8\oplus\Z/8$, depicted in
Table \ref{chart}, but in this case the total obstruction is entirely $KO$-theoretic, as described in
\cite[Prop 4.6]{Cr1}. See also row 3 of \cite[Table 4.1]{Cr1}, which states explicitly that the kernel
of reduction from total obstruction to $KO$-theoretic  obstruction is zero. From our viewpoint,
the non-$KO$-theoretic obstructions are irregular classes such as the one in position $(16p,2)$ in Diagram
\ref{chart}, which must be dealt with in the proof of Proposition \ref{E5}.

Thus $CP^{8p+2}$ lifts to $\bspt(16p-8)$. Since, by Diagram \ref{chart}, $\pi_{16p+6}(\Sigma P_{16p-8})=0$,
this lifting extends over $CP^{8p+3}$, as required for Theorem \ref{mainthm}.

We complete the proof of Theorem \ref{mainthm} by proving Proposition \ref{E5}. We will use the
$bo$-primary classifying spaces $\boo(m)$ constructed in \cite{DM1}. There is a map of fibrations through dimension
$2m-2$
$$\begin{CD}P_m@>>> P_m\w bo\\
@VVV @VVV\\
\bspt(m)@>>> \boo(m)\\
@VVV @VVV\\
\bsp@>=>>\bsp,\end{CD}$$
and there are natural maps of MPTs for these fibrations. We will consider the maps of MPTs
for the following spaces over $\bsp$.
\begin{equation}\label{4spaces}\begin{CD}\bspt(16p-11)@>>>\bspt(16p-9)\\@VVV @VVV\\ \boo(16p-11)@>>> \boo(16p-9).\end{CD}
\end{equation}
We depict in Diagram \ref{chart2} the portion of the Adams spectral sequences in dimensions
$\equiv0$ mod 4 (which is all that is relevant for maps from $HP^{4p+1}$) for $\Sigma P_{16p-11}$,
$\Sigma P_{16p-9}$, $\Sigma P_{16p-11}\w bo$, and $\Sigma P_{16p-9}\w bo$, which correspond to the
$k$-invariants for liftings to each of the spaces in (\ref{4spaces}).

\begin{center}
\begin{minipage}{6.5in}
\begin{diag}\label{chart2}{\bf Possible obstructions for liftings}
\begin{center}

\begin{picture}(350,220)
\def\mp{\multiput}
\def\elt{\circle*{3}}
\mp(5,10)(200,0){2}{\line(1,0){130}}
\mp(5,130)(200,0){2}{\line(1,0){130}}
\mp(10,10)(0,10){3}{\elt}
\mp(10,130)(0,10){3}{\elt}
\mp(10,10)(0,120){2}{\line(0,1){20}}
\mp(50,10)(0,10){4}{\elt}
\mp(50,140)(0,10){3}{\elt}
\mp(50,10)(200,0){2}{\line(0,1){30}}
\mp(50,140)(200,0){2}{\line(0,1){20}}
\mp(250,10)(0,10){4}{\elt}
\mp(250,140)(0,10){3}{\elt}
\mp(210,10)(0,120){2}{\elt}
\mp(90,10)(0,10){7}{\elt}
\mp(90,130)(0,10){7}{\elt}
\mp(90,10)(0,120){2}{\line(0,1){60}}
\mp(290,10)(0,10){5}{\elt}
\mp(290,130)(0,10){5}{\elt}
\mp(290,10)(0,120){2}{\line(0,1){40}}
\put(293,150){\elt}
\mp(330,10)(0,10){8}{\elt}
\put(330,10){\line(0,1){70}}
\mp(330,180)(0,10){3}{\elt}
\put(330,180){\line(0,1){20}}
\mp(130,10)(0,10){8}{\elt}
\put(130,10){\line(0,1){70}}
\put(130,150){\elt}
\put(130,170){\elt}
\mp(130,180)(0,10){3}{\elt}
\put(130,180){\line(0,1){20}}
\mp(-10,0)(0,120){2}{$16p-8$}
\mp(190,0)(0,120){2}{$16p-8$}
\mp(40,0)(0,120){2}{$-4$}
\mp(240,0)(0,120){2}{$-4$}
\mp(88,0)(0,120){2}{$0$}
\mp(288,0)(0,120){2}{$0$}
\mp(128,0)(0,120){2}{$4$}
\mp(328,0)(0,120){2}{$4$}
\put(0,60){$\Sigma P_{16p-11}\w bo$}
\put(200,60){$\Sigma P_{16p-9}\w bo$}
\put(0,180){$\Sigma P_{16p-11}$}
\put(200,180){$\Sigma P_{16p-9}$}
\mp(160,50)(0,120){2}{$\to$}
\mp(65,93)(200,0){2}{$\downarrow$}
\end{picture}
\end{center}
\end{diag}
\end{minipage}
\end{center}

\medskip
By \cite[1.8]{DM1}, $nH_t^{Sp}$ lifts to $\boo(m)$ if and only if for all $i\le t$,
$\nu(\binom ni)\ge \nu(\pi_{4i}(\Sigma P_m\w bo))$. 
By standard methods, one finds
$$\nu(\tbinom{-(4p+2)}{4p+\eps})=\begin{cases}\a(p)-1&\eps=-2\\
2+\a(p)+\nu(p)&\eps=-1\\
\a(p)&\eps=0\\
\a(p)+1&\eps=1.\end{cases}$$
We have $\a(p)=5$. Thus $-(4p+2)H_{4p-1}^{Sp}$ lifts to $\boo(16p-11)$,
and $-(4p+2)H_{4p}^{Sp}$ lifts to $\boo(16p-9)$. By considering the induced map of MPTs
for $\boo(16p-11)\to\boo(16p-9)$, we deduce that $-(4p+2)H_{4p}$ lifts to $E_5$ in the MPT
for $\boo(16p-11)$. Then, since all $k$-invariants for $\bspt(16p-11)$ which are relevant
for $HP^{4p}$ map injectively to those of $\boo(16p-11)$, we infer that $HP^{4p}$ lifts to
$E_5$ of the MPT for $\bspt(16p-11)$. We follow this into $E_5$ of the MPT for $\bspt(16p-9)$.
Since this MPT has no $k$-invariants in dimension $16p+4$ in filtration less than 5, the map
$HP^{4p}\to E_5(16p-9)$ extends over $HP^{4p+1}$, establishing Proposition \ref{E5}.

\section{Proof of Theorem \ref{thm2}}\label{thm2sec}

\begin{pf*}{Proof of Theorem \ref{thm2}} Let $n=2\ell$ with $\a(\ell)=3$. We must show that the map $CP^{2\ell}@>f>> BU@>>> BSO$
which classifies the stable normal bundle $-(2\ell+1)H^\C_{2\ell}$ factors through $BSO(4\ell-4)$.
The fiber in $P_{4\ell-4}\to BSO(4\ell-4)\to BSO$ has an ASS chart that looks like the first four dimensions
of Diagram \ref{chart}, except that if $\ell$ is even, 2 times the bottom class in what appears in that chart as
dimension $16p-4$ equals the sum of the two dots in the box above it. 
We will show that the map lifts to level 3 in the MPT for this fibration, and that the level-3 $k$-invariant is in primary indeterminacy, which implies that the lifting exists. The reason that we did not notice this result in \cite{DM} is apparently
that we were hesitant to consider liftings to $BSO(m)$ when $m$ is divisible by 4 and the bundle is an odd
multiple of the complex Hopf bundle. (See Tables 1.8 and 1.9 of \cite{DM}.)

We let $\but(m)$ denote the classifying space for stably almost complex vector bundles of real geometric
dimension $m$; i.e., it is the pullback of $BU$ and $BO(m)$ over $BO$. As in \cite{DM}, we use spaces
$\but(m)\to\boo_m\to \buo_m$ over $BU$ with fibers through dimension $2m-2$ given by $P_m\to P_m\w bo\to P_m\w bu$.
We need the following charts of homotopy groups.
\begin{center}
\begin{minipage}{6.5in}
\begin{diag}\label{chart3}{\bf Some homotopy groups}
\begin{center}

\begin{picture}(190,180)
\def\mp{\multiput}
\def\elt{\circle*{3}}
\mp(0,15)(0,90){2}{\line(1,0){80}}
\mp(130,15)(0,90){2}{\line(1,0){60}}
\mp(0,0)(0,90){2}{$-4$}
\mp(68,0)(0,90){2}{$0$}
\mp(140,0)(0,90){2}{$-2$}
\mp(178,0)(0,90){2}{$0$}
\mp(10,15)(0,90){2}{\elt}
\mp(40,15)(0,15){2}{\elt}
\mp(40,15)(140,0){2}{\line(0,1){15}}
\mp(180,15)(0,15){2}{\elt}
\mp(70,15)(0,15){3}{\elt}
\put(70,15){\line(0,1){30}}
\put(150,15){\elt}
\mp(40,120)(15,15){3}{\elt}
\put(40,120){\line(1,1){30}}
\mp(70,105)(0,15){3}{\elt}
\put(70,105){\line(0,1){45}}
\mp(150,105)(15,15){3}{\elt}
\put(150,105){\line(1,1){30}}
\mp(180,105)(0,15){2}{\elt}
\put(180,105){\line(0,1){30}}
\mp(100,30)(0,90){2}{$\to$}
\mp(55,75)(105,0){2}{$\downarrow$}
\put(0,55){$\pi_{4\ell+x}(\Sigma P_{4\ell-5}\w bu)$}
\put(130,55){$\pi_{4\ell+x}(\Sigma P_{4\ell-3}\w bu)$}
\put(0,160){$\pi_{4\ell+x}(\Sigma P_{4\ell-5}\w bo)$}
\put(130,160){$\pi_{4\ell+x}(\Sigma P_{4\ell-3}\w bo)$}
\end{picture}
\end{center}
\end{diag}
\end{minipage}
\end{center}
\medskip
We also need the easy fact that, for $\eps=0,1,2$, $\nu(\binom{-(2\ell+1)}{2\ell-\eps})=3,2,2$, respectively.
We use \cite[1.7b]{DM}, which states that, if $p$ is odd, then $pH_n^\C$ lifts to $\boo_m$ if and only if
for all $i\le n$, $\nu(\binom pi)\ge\nu(\pi_{2i}(\Sigma P_m\w bu))$, and, for all even $i\le n$,
$\nu(\binom pi)\ge\nu(\pi_{2i}(\Sigma P_m\w bo))$. This implies that our map $f:CP^{2\ell}\to BU$ lifts to $\boo_{4\ell-3}$, and
$f|CP^{2\ell-1}$ lifts to $\boo_{4\ell-5}$. Thus $f$ lifts to level 3 in the MPT for $\boo_{4\ell-5}$.
By \cite[Tables 8.4,8.12]{Mem}, for $*\le 4\ell$, $\pi_*(\Sigma P_{4\ell-5})\to\pi_*(\Sigma P_{4\ell-5}\w bo)$
is surjective with kernel consisting of a single class in $*=4\ell-1$. Since $H^{4\ell-1}(CP^{2\ell})=0$, we conclude that $f$ lifts to level 3 in the MPT
for $\but(4\ell-5)$, and hence also in those for $\but(4\ell-4)$ and $BSO(4\ell-4)$, using the maps of MPTs induced by
$\but(4\ell-5)\to \but(4\ell-4)\to BSO(4\ell-4)$.

We now employ a standard indeterminacy argument, as explained clearly in \cite{Sin},
to show that the final obstruction in $H^{4\ell}(E_3)$ can be varied, if necessary. Let $E_i$ denote the
spaces in the MPT of $BSO(4\ell-4)\to BSO$. The fiber $F$ of $E_3\to E_2$ is $K_{4\ell-4}\times K_{4\ell-2}\times
K_{4\ell-1}\times K_{4\ell-1}$, corresponding to elements at height 2 in Diagram \ref{chart}, desuspended once.
Here $K_i=K(\Z/2,i)$. If $f_3:CP^{2\ell}\to E_3$ sends the $k$-invariant $k_3\in H^{4\ell}(E_3)$ trivially, then the lifting to $BSO(4\ell-4)$ exists, since there are no more even-dimensional $k$-invariants. If $f_3^*(k_3)\ne 0$, then we will show that the composite
\begin{equation}\label{comp}CP^{2\ell}@>\iota_{4\ell-4}\times f_3>>F\times E_3@>\mu>>E_3,\end{equation}
which is also a lifting of $f$,
sends $k_3$ to 0, and hence the lifting to $BSO(4\ell-4)$ exists. Here $\mu$ denotes the action of the fiber on the total space
in the principal fibration, and $\iota_{4\ell-4}:CP^{2\ell}\to F$ is the map which is nontrivial into the first factor
of $F$. This will follow because a computation of the relations in the MPT, performed below, shows that
\begin{eqnarray}\mu^*(k_3)&=&1\times k_3+\sq^1\io_{4\ell-1}'\times 1+\sq^2\io_{4\ell-2}\times1+\io_{4\ell-2}\times w_2\nonumber\\
&&+\sq^4\io_{4\ell-4}\times1+\io_{4\ell-4}\times(w_4+w_2^2).\label{MPTreln}\end{eqnarray}
Thus the composite (\ref{comp}) sends $k_3$ to $$f_3^*(k_3)+\sq^4(x^{2\ell-2})+x^{2\ell-2}\cdot(w_4+w_2^2)(-(2\ell+1)H).$$
Here $x$ denotes the generator of $H^2(CP^{2\ell})$. Since $w_2(-(2\ell+1))H=x$, and $\sq^4(x^{2\ell-2})$ and $w_4(-(2\ell+1))$ are either both nonzero ($\ell$ even) or both zero ($\ell$ odd), we deduce that $f_3^*(k_3)$ can be varied, if necessary, establishing the lifting.

We conclude the proof by listing the relations in the MPT of $BSO(4\ell-4)\to BSO$, the last of which yields
the crucial fact (\ref{MPTreln}). These are computed by the method initiated in \cite{GM} and utilized in such 
papers as \cite{DM}, \cite{Ran}, and \cite{Sin}. It is a matter of building a minimal resolution using Massey-Peterson
algebras. In this table, $\eps=1$ if $\ell$ is even, and 0 if $\ell$ is odd. 

\begin{center}
\renewcommand{\arraystretch}{1.2}
\begin{tabular}{|ll|}
\hline
$w_{4\ell-3}$&\\
$w_{4\ell-2}$&\\
$w_{4\ell}$&\\
\hline
$k^1_{4\ell-3}$:&$\sq^1w_{4\ell-3}$\\
$k^1_{4\ell-2}$:&$(\sq^2+w_2)w_{4\ell-3}$\\
$k^1_{4\ell-1}$:&$(\sq^2+w_2)w_{4\ell-2}$\\
$k^1_{4\ell}$:&$\sq^1w_{4\ell}+(\sq^2+w_2)\sq^1w_{4\ell-2}$\\
$\kt^1_{4\ell}$:&$\eps\sq^1w_{4\ell}+(\sq^4+w_4)w_{4\ell-3}+(w_2\sq^1+w_3)w_{4\ell-2}$\\
\hline
$k^2_{4\ell-3}$:&$\sq^1k^1_{4\ell-3}$\\
$k^2_{4\ell-1}$:&$(\sq^2+w_2)k^1_{4\ell-2}+(\sq^3+w_3)k^1_{4\ell-3}$\\
$k^2_{4\ell}$:&$\sq^1k^1_{4\ell}+(\sq^2+w_2)k^1_{4\ell-1}$\\
$\kt^2_{4\ell}$:&$\sq^1\kt^1_{4\ell}+(\sq^2\sq^1+w_3)k^1_{4\ell-2}+(\sq^4+w_4+w_2^2+w_2\sq^2)k^1_{4\ell-3}$\\
\hline
$k^3_{4\ell-3}$:&$\sq^1k^2_{4\ell-3}$\\
$k^3_{4\ell}$:&$\sq^1\kt^2_{4\ell}+(\sq^2+w_2)k^2_{4\ell-1}+(\sq^4+w_4+w_2^2)k^2_{4\ell-3}$\\
\hline
\end{tabular}
\end{center}
\end{pf*}

\section{Discussion of one of Crabb's proofs}\label{sec4}

In \cite[0.2]{Crabb}, Crabb presented many new immersions of complex projective spaces.
About his proof, he wrote \lq\lq Details will be omitted," although sketched arguments for
each case were presented. We have checked the details of
his arguments, and found what appears to be a flaw in one case. This case was of particular
interest to us, because, if true, it would have implied a new immersion result for real
projective space which would be an addition to \cite{immtable}.
We present here our analysis of this case.
\begin{prop} The argument for the portion of \cite[0.2]{Crabb} which states that if $k$ is even and $\a(k)=4$ then $CP^{2k+1}$\label{unresolved}
immerses in $\R^{8k-6}$ is invalid. The question of whether this immersion exists is unresolved.\end{prop}
\begin{pf} Let $k=2K$. To obtain the immersion, one must prove that $-(4K+2)H_{4K+1}^\C$ lifts to $\bspt(8K-8)$. The $KO$-theoretic
obstruction for lifting this bundle, calculated similarly to the one in Section \ref{intro}, is 0.
However, there are several elements in the kernel of the reduction from the total obstruction to the
$KO$-theoretic obstruction which cannot be ruled out. In \cite[Table 4.1]{Cr1}, these are the
three $\Z/2$'s in row 1, column 9. In our Diagram \ref{chart}, which, after a slight reindexing, serves
as the obstructions for this lifting question as well as the one considered in Section \ref{intro}, these
correspond to three of the many dots in column $16p+2$.

Crabb realized that these could cause a problem, and so he hoped to utilize the factorization through $HP^{2K}$,
similarly to what we did in Section \ref{intro}. In fact, he wrote in his proof on page 166 that in several
cases, including this one, it can be shown that the bundle over the quaternionic projective
space is stably equivalent to a bundle of the desired dimension. In this case, he would be saying that
$-(2K+1)H_{2K}^{Sp}$ lifts to $\bspt(8K-8)$. However, this lifting does not exist; its $KO$-theoretic obstruction is
nonzero. 

This can be seen similarly to our calculation in Section \ref{intro}. We use \cite[Prop 2.5]{Crabb}. With $\lambda(T)$ as in Section \ref{intro}, the necessary condition is
that there exists an integer $e$ such that $e(\l^{4K-1})_2\equiv 16$ mod 32, $e(\l^{4K-1})_1\equiv8$ mod 16, and
$e(\l^{4K-1})_0\equiv0$ mod 2. Since
$$(1-\tfrac13T+\tfrac 8{45}T^2)^{4K-1}=1+\text{od}\cdot T+\text{od}\cdot T^2$$
mod $(T^3)$, where od denotes an odd fraction, we require $e$ to satisfy both $e\equiv16$ mod 32, and $e\equiv 8$ mod 16, which is clearly impossible.
\end{pf}

\smallskip
We close by commenting on the relationship between Crabb's necessary conditions for
immersion involving powers of $\l(T)=(\sinh^{-1}(\sqrt T)/\sqrt T)^2$, which we have used above, and the Sigrist-Suter
necessary condition involving powers of $\log(1+t)/t$ in Theorem \ref{SigSthm}. These conditions involving power series
can be directly related to one another by a slight extension of \cite[1.5]{BD}, which was proved using a proof suggested to the authors by Crabb, and based on his
earlier topological work described in \cite[\S4]{CK}. 

\section{Evaluation of a Sigrist-Suter condition}\label{sec5}
In Theorem \ref{SigSthm}, a general statement of a necessary condition for $CP^n$ to immerse in $\R^{4n-2\a(n)}$
is presented. We evaluate this explicitly when $\a(n)=8$ in the following result, which we have incorporated into Tables \ref{evtabl} and \ref{odtabl}.
\begin{prop} If $\a(n)=8$ and $CP^n$ immerses in $\R^{4n-2\a(n)}$, then $\nu(n-7)=3$ or $\nu(n-6)=4$ or $\nu(n-5)=5$ or $\nu(n-4)=7$.\end{prop}
This follows readily from Theorem \ref{SigSthm} and the following lemma.
\begin{lem} Let $(\log(1+t)/t)^m=\sum c_{m,i}t^i$. Then
$$\nu(c_{m,7})=\nu(c_{m,8})<\nu(c_{m,i})\text{ for all }i<7$$
if and only if $\nu(m-7)=4$ or $\nu(m-5)=5$ or $\nu(m-3)=6$ or $\nu(m-1)=8$.\end{lem}
\begin{pf} Let
$$v(m)=(v_0(m),\ldots,v_1(m))=(\nu(c_{m,0}),\ldots,\nu(c_{m,8})).$$
Also define a function, for which some of its values are only specified to satisfy an inequality, by
$$\nu(k,e)=\begin{cases}\nu(k)&\text{if }\nu(k)<e\\
\ge e&\text{if }\nu(k)\ge e.\end{cases}$$
The lemma follows immediately from the following, which we will prove for $v(m)$ by induction on $m$.
Note that some components are only asserted to satisfy an inequality.
\begin{enumerate}
\item If $e\ge3$, then $v(2^e(2a+1))=(0,e-1,e-3,e-3,e-6,e-5,e-7,e-7,e-11)$.
\item $v(7+8k)=(0,-1,-2,-3,-4,-5,-6,-7,-8+\nu(k,2))$.
\item $v(5+8k)=(0,-1,-1,\ge-1,-4,-5,\ge-4,-6,-8+\nu(k,3))$.
\item $v(3+32k)=(0,-1,-2,-3,\ge-1,-4,-4,-5,-6+\nu(k,2))$.
\item $v(19+32k)=(0,-1,-2,-3,-2,-4,\ge-3,\ge-4,-7)$.
\item $v(11+16k)=(0,-1,-2,-3,-3,\ge-3,-5,-6,-8)$.
\item $v(9+16k)=(0,-1,\ge-1,-2,-3,-4,-4,\ge-3,-8)$.
\item $v(17+32k)=(0,-1,0,-2,-2,-3,-3,\ge-2,-7)$.
\item $v(1+32k)=(0,-1,0,-2,\ge-1,\ge-2,\ge-2,-3,-6+\nu(4,k))$.
\end{enumerate}

We begin by using {\tt Maple} to verify $v(m)$ for $m=8$, 7, 5, 3, 19, 11, 9, 17, and 1. 
The induction proof for $v(2^e)$ is obtained using
$$\sum c_{2^{e+1},i}t^i=(\sum c_{2^e,i}t^i)^2.$$
We have
\begin{equation}\label{veq}v_i(2^{e+1})=1+v_i(2^e),\end{equation}
since, as is easily verified, the RHS of (\ref{veq}) is strictly less than $1+v_j(2^e)+v_{i-j}(2^e)$ for
$1\le j\le i/2$.

Next we obtain the claim for $v(2^e(2a+1))$ from $(\sum c_{2^e,i}t^i)(\sum c_{2^{e+1}a,i}t^i)$.
From this, we obtain $v_i(2^e(2a+1))=v_i(2^e)$ since the strict minimum of $v_j(2^e)+v_{i-j}(2^{e+1}a)$ is obtained
when $j=i$, essentially the same verification as the previous one.

The cases in which the asserted $v_i(u+2^ek)$ is not of the form $\ge t$ or $\nu(k,t)$ are, except for $v_7(1+32k)$ with $k$ odd, obtained from
$v_i(u)+v_0(2^ek)$ which, in these cases, is strictly less than other values of $v_j(u)+v_{i-j}(2^ek)$.
Quite a few verifications are required for this. For example, since $v_4(2^ek)=e+\nu(k)-6$, it is relevant that
$$v_i(u)<v_{i-4}(u)-\begin{cases}3&u=7,5\\
2&u=11,9\\
1&u=3,19,17,1\end{cases}$$
in the cases in which the asserted value of $v_i(u)$ is a single integer. These cases are differentiated
according to whether the argument of $v$ is $u+8k$, $u+16k$, or $u+32k$. If $k$ is odd, then $v_7(1+32k)=-3$
comes from having each of $v_7(1)+v_0(32k)$, $v_3(1)+v_4(32k)$, and $v_1(1)+v_6(32k)$ equal $-3$.

In the cases in which the asserted $v_i(u+2^ek)$ is of the form $\ge t$, one verifies, unless $i=7$ and $u=9$ or 17, that $v_j(u)+v_{i-j}(2^ek)\ge t$ for all $j$, with the possibility of equality for several values of $j$.
For example, $v_6(5+8k)\ge -4$ comes from $v_6(5)+v_0(8k)\ge-4$, $v_4(5)+v_2(8k)\ge-4$, and $v_2(5)+v_4(8k)\ge-4$,
while $v_j(5)+v_{6-j}(8k)>-4$ for $j=0$, 1, 3, and 5. If $k$ is odd, the exceptional case $v_7(9+16k)\ge -3$ comes from
$v_1(9)+v_6(16k)=-4$ and $v_3(9)+v_4(16k)=-4$, with other values $\ge-3$. The case $v_7(17+32k)$ is similar.

Finally, for $\nu(k)=0$ or 1, $v_8(7+8k)=-8+\nu(k)$ comes from $v_0(7)+v_8(8k)$, which is strictly less than
all other $v_j(7)+v_{8-j}(8k)$, while if $\nu(k)\ge 2$, we have $v_0(7)+v_8(8k)=\nu(k)-8$ and $v_8(7)+v_0(8k)\ge-6$
with other terms larger. The same argument works for $v_8(5+8k)$, $v_8(3+32k)$, and $v_8(1+32k)$.
\end{pf}

\def\line{\rule{.6in}{.6pt}}

\end{document}